\documentclass[12pt]{article}

\usepackage{amssymb}
\usepackage{amsthm}

\newtheorem*{theorem}{Theorem}
\newtheorem*{definition}{Definition}

\newtheorem{lemma}{Lemma}

\newenvironment{remark}[1][Remark]
           {\medbreak\noindent \textbf{#1 \enspace}}
           {\par \medbreak}

\newcommand{\Card}{\mathop{\rm Card}\nolimits}
\newcommand{\ORD}{\mathop{\rm ORD}\nolimits}

\newcommand{\A}{{\mathcal A}}
\newcommand{\QQ}{{\mathbb Q}}

\begin{document}
\baselineskip=18pt

\thispagestyle{empty}

\title{David's Trick}

\author{Sy D. Friedman\thanks{%
  Research supported by NSF Contract \#9625997-DMS}\\[1ex]
  {\normalsize M.I.T.}
}
\date{}

\maketitle

\thispagestyle{empty}

In David~\cite{82} a method is introduced for creating reals $R$
which not only code classes in the sense of Jensen coding but in
addition have the property that in $L [R]$, $R$ is the unique
solution to a $\Pi^1_2$ formula.  In this article we cast David's
``trick'' in a general form and describe some of its uses.

\begin{theorem}
  Suppose $A \subseteq \ORD$, $\langle L [A], A \rangle \models Z
  F C + 0^\#$ does not exist and suppose that for every infinite
  cardinal $\kappa$ of $L [A]$, $H_\kappa^{L [A]} = L_\kappa [A]$
  and $\langle L_\kappa [A]$, $A \cap \kappa \rangle \models
  \varphi$.  Then there exists a $\Pi^1_2$ formula $\psi$ such
  that:

  \renewcommand{\theenumi}{\alph{enumi}}%
  \renewcommand{\labelenumi}{(\theenumi)}%
  \begin{enumerate}
  \item 
    If $R$ is a real satisfying $\psi$ then there is $A \subseteq
    \ORD$ as above, definable over $L [R]$ in the parameter $R$.
   
  \item 
    For some tame, $\langle L [A], A \rangle$-definable,
    cofinality-preserving forcing $P$, $P \Vdash \exists R
    \psi (R)$.
  \end{enumerate}
  
  Moreover if $A$ preserves indiscernibles then $\psi$ has a
  solution in $L [A, 0^\#]$, preserving indiscernibles.

\end{theorem}

\begin{remark}
  \renewcommand{\theenumi}{\arabic{enumi}}%
  \renewcommand{\labelenumi}{(\theenumi)}%
  \begin{enumerate}
  \item 
    We require that $H_\kappa^{L [A]}$ equal $L_\kappa [A]$ for
    infinite $L [A]$-cardinals solely to permit
    cofinality-preservation for $P$; if cofinality-preservation
    is dropped then such a requirement is unnecessary, by coding
    $A$ into $A^*$ with this requirement and then applying our
    result to $A^*$.

  \item 
    A class $A$ \emph{preserves indiscernibles} if the Silver
    indiscernibles are indiscernible for $\langle L [A], A
    \rangle$.  It follows from the technique of Theorem~0.2 of
    Beller-Jensen-Welch~\cite{82} (see Friedman~\cite{98}) that if $A$
    preserves indiscernibles then $A$ is definable from a real $R
    \in L [A, 0^\#]$, preserving indiscernibles.
  \end{enumerate}
\end{remark}

\begin{proof}
  Our plan is to create an $\langle L [A], A \rangle$-definable,
  tame, cofinality-preserving forcing $P$ for adding a real $R$
  such that whenever $L_\alpha [R] \models Z F^{-}$ there is
  $A_\alpha \subseteq \alpha$, definable over $L_\alpha [R]$ (via
  a definition independent of $\alpha$) such that $L_\alpha [R]
  \models$ for every infinite cardinal $\kappa$, $H_\kappa =
  L_\kappa [A_\kappa]$ and $\varphi$ is true in $\langle L_\kappa
  [A_\alpha], A_\alpha \cap \kappa \rangle$.  This property
  $\psi$ of $R$ is $\Pi^1_2$ and gives us (a), (b) of the
  Theorem.  The last statement of the Theorem will follow using
  Remark (2) above.

$P$ is obtained as a modification of the forcing from
Friedman~\cite{97}, used to prove Jensen's Coding Theorem (in the
case where $0^\#$ does not exist in the ground model).  The
following definitions take place inside $L [A]$.

\begin{definition}[Strings]
  Let $\alpha$ belong to $\Card =$ the class of all infinite
  cardinals.  $S_\alpha$ consists of all $s: [\alpha, \left| s
  \right| ) \to 2$, $\alpha \le \left| s \right| < \alpha^{+}$
  such that $\left| s \right|$ is a multiple of $\alpha$ and:

  \renewcommand{\theenumi}{\alph{enumi}}%
  \renewcommand{\labelenumi}{(\theenumi)}%
  \begin{enumerate}
  \item 
    $\eta \le \left| s \right| \to L_\delta [A \cap \alpha, s
    \upharpoonright \eta ] \models \Card \eta \le \alpha$ for
    some $\delta < (\eta^{+})^L \cup \omega_2$.

  \item 
    If $\A = \langle L_\beta [ A \cap \alpha, s \upharpoonright
    \eta ], s \upharpoonright \eta \rangle \models (Z F^{-}
    \hbox{ and } \eta = \alpha^{+})$ then over $\A$, $s
    \upharpoonright \eta$ codes a predicate $A (s \upharpoonright
    \eta, \beta) = A^* \subseteq \beta$ such that $A^* \cap
    \alpha = A \cap \alpha$ and for every cardinal $\kappa$ of
    $L_\beta [A^*]$, $H^{L_\beta [A^*]}_\kappa = L_\kappa [A^*]$
    and $\langle L_\kappa [A^*], A^* \cap \kappa \rangle \models
    \varphi$. 
  \end{enumerate}
\end{definition}

\begin{remark}
  When in (b) above we say that $s \upharpoonright \eta$
  \emph{codes} $A^*$ we are referring to the canonical coding
  from the proof of Theorem~4 of Friedman~\cite{97} of a subset of
  $\beta$ by a subset of $(\alpha^{+})^\A = \eta$ (relative to
  $A \cap \alpha$).

  The remainder of the definitions from the proof of Theorem~4 of
  Friedman~\cite{97} remain the same in the present context.  We
  now verify that he proofs of the lemmas from Friedman~\cite{97}
  can successfully accommodate the new restriction (clause (b))
  on elements of $S_\alpha$.
\end{remark}

\begin{lemma}[Distributivity for $R^s$]
  \label{lem:1}
  Suppose $\alpha \in \Card$,
  $s \in S_{\alpha^{+}}$.  Then $R^s$ is
  $\alpha^{+}$-distributive in $\A^s$.
\end{lemma}

\begin{proof}
  Proceed as in the proof of Lemma~5 of Friedman~\cite{97}.  The
  only new point is to verify that in the proof of the Claim,
  $t_\lambda$ satisfies clause (b) (of the new definition of
  $S_\alpha$).  The fact that $s$ belongs to $S_{\alpha^{+}}$ and that
  $t_\lambda$ codes $\bar{H}_\lambda$ imply that clause (b) holds
  for $t_\lambda$ whenever $\beta$ is at most $\bar{\mu}_\lambda
  = $ the height of $\bar{H}_\lambda$.  But as $\left| t_\lambda
  \right|$ is definably singular over $L_{\bar{\mu}_\lambda}
  [t_\lambda]$ these are the only $\beta$'s that concern~us.
\end{proof}

\begin{lemma}[Extendibility of $P^s$]
  \label{lem:2}
  Suppose $p \in P^s$, $s \in S_\alpha$, $X \subseteq \alpha$, \hbox{$X
  \in \A^s$}.  Then there exists $q \le p$ such that $X \cap \beta
  \in \A^{q_\beta}$ for each $\beta \in \Card \cap \alpha$.
\end{lemma}

\begin{proof}
  Proceed as in the proof of Lemma~6 of Friedman~\cite{97}.  In the
  definition of $q$, the only instances of clause (b) to check
  are for $s_\beta$ when Even $(Y \cap \beta)$ codes $s_\beta$, $
 s_\beta$ satisfying clause (a) of the definition of membership
 in $S_\beta$.  But the embedding $\bar{\A}_\beta \to \A$ is
 $\Sigma_1$-elementary and instances of clause (b) refer to
 ordinals less than the height of $\A$; so the fact that $s$
 belongs to $S_\alpha$ implies that $s_\beta$ belongs to
 $S_\beta$. 
\end{proof}

\begin{lemma}[Distributivity for $P^s$]
  \label{lem:3}
  Suppose $s \in S_{\beta^{+}}$, $\beta \in \Card$.

  \renewcommand{\theenumi}{\alph{enumi}}%
  \renewcommand{\labelenumi}{(\theenumi)}%
  \begin{enumerate}
  \item 
    If $\langle D_i \mid i < \beta \rangle \in \A^s, D_i$ $i^+$
    dense on $P^s$ for each $i < \beta$ and $p \in P^s$ then
    there is $q \le p$, $q$ meets each $D_i$. 

  \item 
    If $p \in P^s$, $f$ small in $\A^s$ then there exists $q \le
    p$, $q \in \Sigma^p_f$.
  \end{enumerate}
\end{lemma}

\begin{proof}
  Proceed as in the proof of Lemma~7 of Friedman~\cite{97}.  In the
  Claim we must verify that $p^\lambda_\gamma$ satisfies
  clause~(b).  But once again this is clear by the
  $\Sigma_1$-elementary of $\bar{H}_\lambda (\gamma)$ and the
n  fact that $L_{\bar{\mu}} [A \cap \gamma, p^\lambda_\gamma]
  \models \left| p^\lambda_\gamma \right|$ is
  $\Sigma_1$-singular, where $\bar{\mu} = $ height of
  $\bar{H}_\lambda (\gamma)$.
\end{proof}

The argument of the proof of Lemma~\ref{lem:3} can also be
applied to prove the distributivity of $P$, observing that when
building sequences of conditions $\langle p^i \mid i < \lambda
\rangle$, $\lambda$ limit to meet an $\langle L [A], A
\rangle$-definable sequence of dense classes, one has that
$p^\lambda_\gamma$ codes $\bar{H}^\lambda (\gamma)$ of height
$\bar{\mu}$, where $L_{\bar{\mu} + 1} [A \cap \gamma,
p^\lambda_\gamma] \models \left| p^\lambda_\gamma \right|$ is not
a cardinal.  Thus there is no additional instance of clause (b) to
verify beyond those considered in the proof of Lemma 3.

Thus $P$ is tame and cofinality-preserving.  The final statement
of the Theorem also follows, using Remark (2) immediately after
the statement of the Theorem.
\end{proof}

\subsection*{Applications}

  \renewcommand{\theenumi}{\arabic{enumi}}%
  \renewcommand{\labelenumi}{(\theenumi)}%
  \begin{enumerate}
  \item 
    Local $\Pi^1_2$-Singletons.  David~\cite{82} proves the
    following:  There is an $L$-definable forcing $P$ for adding a
    real $R$ such that $R$ is a $\Pi^1_2$-singleton in every
    set-generic extension of $L [R]$ (via a $\Pi^1_2$ formula
    independent of the set-generic extension).  This is
    accomplished as follows:  One can produce an $L$-definable
    sequence $\langle T (\kappa) \mid \kappa \hbox{ an infinite }
    L\hbox{-cardinal} \rangle$ such that $T (\kappa)$ is a
    $\kappa^{++}$-Suslin tree in $L$ for each $\kappa$ and the
    forcing $\prod T (\kappa)$ for adding a branch $b (\kappa)$
    through each $T (\kappa)$ (via product forcing, with Easton
    support) is tame and cofinality-preserving.  Now for each $n$
    let $X_n \subseteq \omega_1^L$ be class-generic over $L$,
    $X_n$ codes a branch through $T (\kappa)$ iff $\kappa$ is of
    the form $(\aleph^L_{\lambda + n})$, $\lambda$ limit.  The
    forcing $\prod P_n$, where $P_n$ adds $X_n$, can be shown to be
    tame and cofinality-preserving.  Finally over $L [\langle X_n
    \mid n \in \omega \rangle]$ add a real $R$ such that $n \in
    R$ iff $R$ codes $X_n$.  Then one has that in $L [R]$, $n \in
    R$ iff $T (\aleph^L_{\lambda + n})$ is not $\aleph^L_{\lambda +
      n}$-Suslin for sufficiently large $\lambda$.  Clearly this
    characterization will still hold in any set-generic extension
    of $L [R]$.  David's trick is used to strengthen this to a
    $\Pi^1_2$ property of $R$.

  \item 
    A Global $\Pi^1_2$-Singleton.  Friedman~\cite{90} produces a
    $\Pi^1_2$-singleton $R$,\break 
    \hbox{$0 <_L R <_L 0^\#$}.  This is
    accomplished as follows:  assume that one has an index for a
    $\Sigma_1 (L)$ classification $(\alpha_1 \cdots \alpha_n)
    \mapsto r (\alpha_1 \cdots \alpha_n)$ that produces $r
    (\alpha_1 \cdots \alpha_n) \in 2^{< \omega}$ for each
    $\alpha_1 < \cdots < \alpha_n$ in ORD such that $R = \cup \{ r
    (i_1 \cdots i_n) \mid i_1 < \cdots < i_n$ in $I =$ Silver
    indiscernibles $\}$.  For each $r \in 2^{< \omega}$ there is
    a forcing $\QQ (r)$ for ``killing'' all $(\alpha_1 \cdots
    \alpha_n)$ such that $r (\alpha_1 \cdots \alpha_n)$ is
    incompatible with $r$.  No $(i_1 \cdots i_n)$ from $I^n$ can be
    killed. Now build $R$ such that $r \subseteq R$ iff $R$ codes a $\QQ (r)$-generic.  Then $R$ is
    the unique real with this property.  David's trick is used to
    strengthen this to a $\Pi^1_2$ property.

  \item 
    New $\Sigma^1_3$ facts.  Friedman~\cite{98} shows that if $M$
    is an inner model of ZFC, $0^\# \notin M$, then there is a
    $\Sigma^1_3$ sentence false in $M$ yet true in a forcing
    extension of~$M$.  This is accomplished as follows:  let
    $\langle C_\alpha | \alpha  \hbox{ $L$-singular} \rangle$ be a
    $\square$-sequence in $L$; i.e., $C_\alpha$ is CUB in
    $\alpha$, $ot C_\alpha < \alpha$, $\bar{\alpha} \in \lim
    C_\alpha \to C_{\bar{\alpha}} = C_\alpha \cap \bar{\alpha}$.
    Define $n (\alpha) = 0$ if $otC_\alpha$ is $L$-regular and
    otherwise $n (\alpha) = n (ot C_\alpha) + 1$.
    Then for some $n$, $\{ \alpha \mid n (\alpha) =n \}$ is
    stationary in $M$.  And for each $n$, there is a tame forcing
    extension of $M$ in which $\{ \alpha \mid n (\alpha) \le n
    \}$ is non-stationary, and is in fact disjoint from the class
    of limit cardinals.  David's trick is used to strengthen the
    latter into a $\Sigma^1_3$ property.
  \end{enumerate}


\frenchspacing
\begin{thebibliography}{99}

\bibitem[82]{82}
  R. David, A Very Absolute $\Pi^1_2$-Singleton, {\it Annals
    of Pure and Applied Logic} {\bf 23} pp.~101-120.

\bibitem[82]{82}
  A. Beller, R. Jensen, P. Welch, {\it Coding the Universe}, book, {\it Cambridge University Press}

\bibitem[90]{90}
  S. Friedman, The $\Pi^1_2$ Singleton Conjecture, {\it Journal of the American Mathematical Society}, Vol.3, No.4, pp. 771-791.

\bibitem[97]{97}
  S. Friedman, Coding without Fine Structure, {\it Journal of Symbolic Logic}, Vol.62, No.3, pp.808-815.

\bibitem[98]{98}
  S. Friedman, New $\Sigma^1_3$ Facts, to appear.

\bibitem[99]{99}
  S. Friedman, {\it Fine Structure and Class Forcing}, book, in preparation.

\end{thebibliography}
\end{document}